\documentclass{article}

\usepackage{arxiv}

\usepackage[utf8]{inputenc} 
\usepackage[T1]{fontenc}    
\usepackage{hyperref}       
\usepackage{url}            
\usepackage{booktabs}       
\usepackage{amsfonts}       
\usepackage{nicefrac}       
\usepackage{microtype}      
\usepackage{lipsum}

\usepackage{amsmath, amsthm, amssymb}    
\usepackage{graphicx}			

\usepackage[percent]{overpic}

\title{A Family of Fern-like Ternary Complex Trees
}

\author{Bernat Espigul\'e
\vspace{10pt}\\
bernat@espigule.com
}

\begin{document}

\maketitle


\begin{abstract}
A ternary complex tree related to the golden ratio is used to show how the theory of complex trees works. We use the topological set of this tree to obtain a parametric family of trees in one complex variable. Even though some real ferns and leaves are reminiscent to elements of our family of study, here we only consider the underlying mathematics. We provide aesthetically appealing examples and a map of the unstable set~$\mathcal{M}$ for this family. Moreover we show that some elements found in the boundary of the unstable set~$\mathcal{M}$ possess interesting algebraic properties, and we explain how to compute the Hausdorff dimension and the shortest path of self-similar sets described by trees found outside the interior of the unstable set~$\mathcal{M}$.
\end{abstract}

\begin{figure}[h!tbp]
	\centering
 \begin{overpic}[width=.876\textwidth,tics=10]{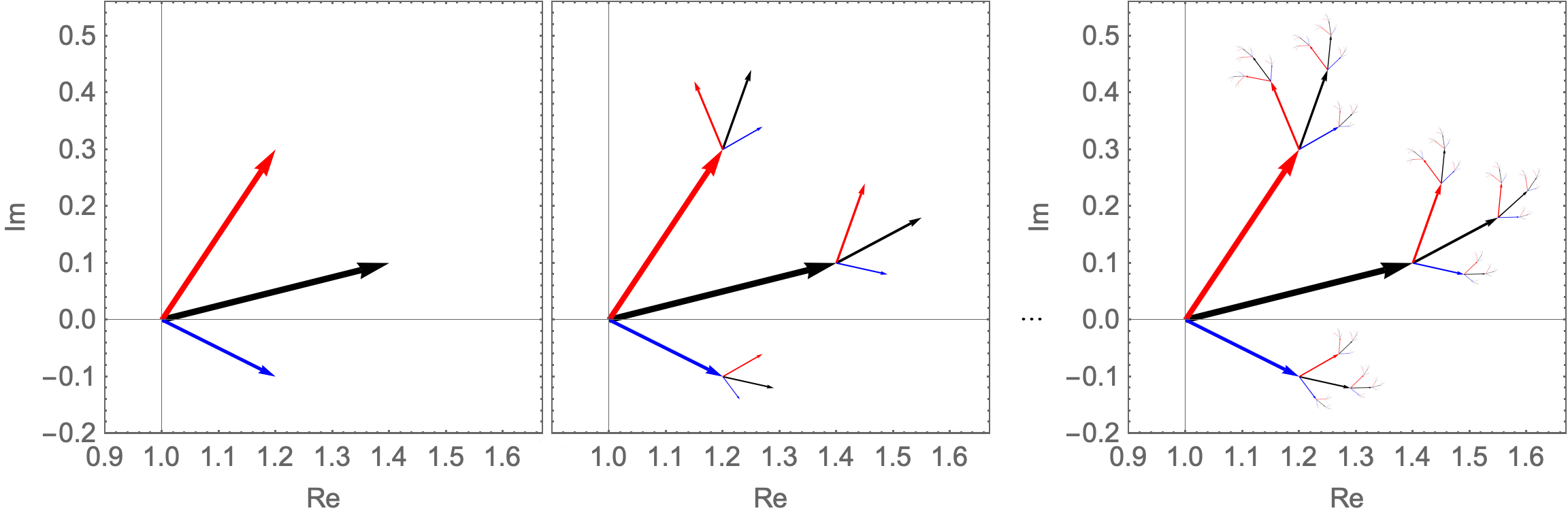}
 \put (7,9) {$\phi(e_0)$}
 \put (25,16) {$\phi(c_1)$}
 \put (16,24) {$\phi(c_2)$}
 \put (18,8) {$\phi(c_3)$}
\end{overpic}
	\caption{Ternary complex tree~$T_A=T{\{c_1,c_2,c_3\}}=T\{0.4+i0.1,0.2+i0.3,0.2-i0.1\}$.}
	\label{fig:0}
\end{figure}
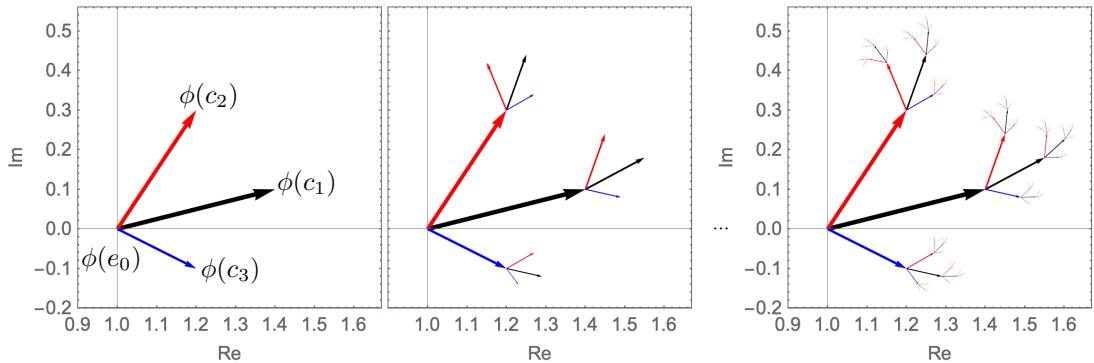

A ternary complex tree~$T_A$ is a fractal tree with all of its branch-nodes encoded by the geometric map~$\phi$ introduced in~\cite{Espigule2019} that sends any word~$w$ composed of~$m$~complex-valued letters $w_1,w_2,\dots,w_m$ taken from a ternary alphabet $A=\{c_1,c_2,c_3\}$ to a complex point~$\phi(w)\in\mathbb{C}$ in the following geometric series style 
\begin{equation*}
	\phi(w):=1+w_1+w_1w_2+\dots+w_1w_2\dots w_m.
\end{equation*}
\noindent If a word $w=w_1w_2w_3\dots w_m$ has a finite sequence of letters, $w_1,w_2,\dots,w_m$, then the point~$\phi(w)$ is called a \textbf{node} of the complex tree~$T_A$. But if a word $w=w_1w_2w_3\dots w_k\dots$ has infinite length, then $\phi(w)$ is called a \textbf{tip point} of the complex tree~$T_A$ and we express it as an infinite sum
 \begin{equation}
 	\phi(w)=1+w_1+w_1w_2+\dots+w_1w_2\dots w_k+\dots=\sum_{k=0}^\infty w_1w_2\dots w_k=\sum_{k=0}^\infty w_{|k}
 \end{equation}
\noindent where each summand $w_{|k}=w_1w_2\dots w_k$ is assumed to be the complex multiplication of the individual letters of the word~$w$ pruned up to its $k$th letter.
For $k=0$, we have that $w_{|0}=w_0=e_0$ where $e_0$ is the \textbf{empty string} with assigned value equal to 1. The node $\phi(e_0)=1$ is called the \textbf{root} of the complex tree where the three first-level nodes, $\phi(c_1)=1+c_1$,~$\phi(c_2)=1+c_2$,~and~$\phi(c_3)=1+c_3$, sprout from, see figure~\ref{fig:0}. 

\begin{figure}[h!tbp]
	\centering
 \begin{overpic}[width=1.08\textwidth,tics=10]{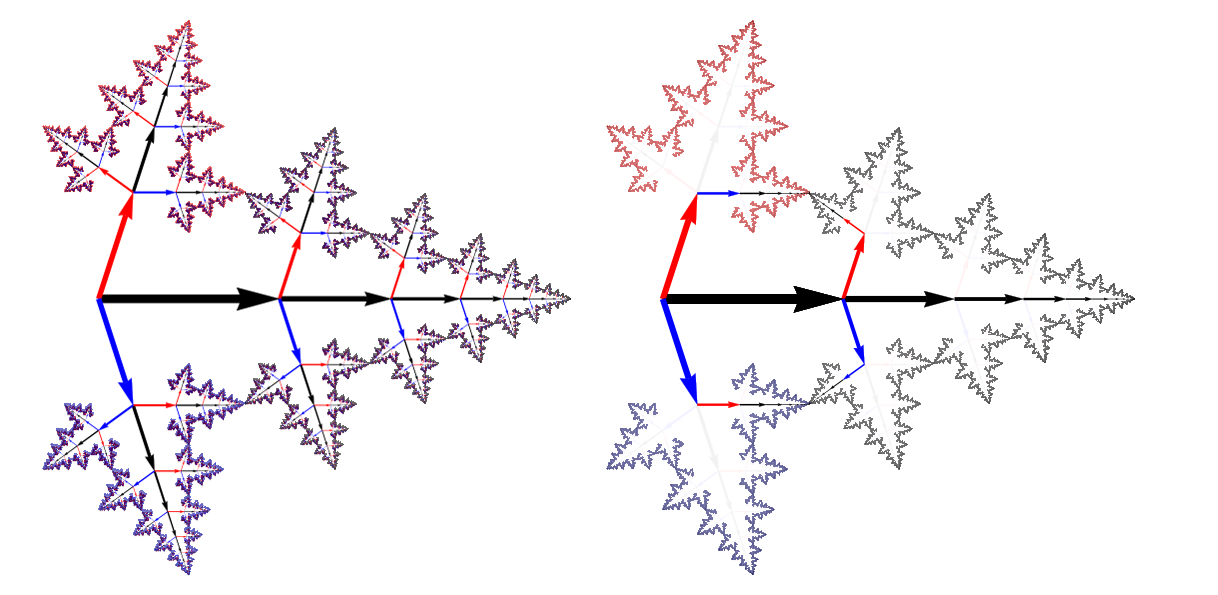}
 \put (94,24) {$\phi(\overline{1})$}
 \put (64,40) {$\phi(23\overline{1})=\phi(122\overline{1})$}
 \put (64,8) {$\phi(32\overline{1})=\phi(133\overline{1})$}
\end{overpic}
	\caption{Complex tree~$T_{\{c_1,c_2,c_3\}}=T\{(-1+\sqrt{5})/2,(-2+\sqrt{5}+i\sqrt{5-2\sqrt{5}})/2,(-2+\sqrt{5}-i\sqrt{5-2\sqrt{5}})/2\}$. The scaling factors for the first-level pieces are $|c_1|=1/\tau$~and~$|c_2|=|c_3|=1/\tau^2$ respectively, where $\tau$ is the golden ratio~$\tau=1.6180\dots$. This tree is mirror-symmetric, $c_1=1/\tau$ and $c_2=c_3^*$.}
	\label{fig:1}
\end{figure}

\noindent By imposing a color code, $\{1\rightarrow c_1, 2\rightarrow c_2,3\rightarrow c_3\}$, a word~$u=u_1u_2\dots u_m$ of a node~$\phi(u)$ can be retrieved by reading the color sequence of the branch-path that goes from the root~$\phi(e_0)$ to the desired node~$\phi(u)$. From now on, letters~$w_k\in\{c_1,c_2,c_3\}$ found in words will be replaced by the numeric symbols~$\{1,2,3\}$ to facilitate reading, for example $c_2c_3c_1c_1=2311$. Infinite sequence of letters $w_1,w_2,\dots,w_k,\ldots\in A$ that are eventually periodic have their associated tip point~$\phi(w)$ reduced to algebraic expressions in terms of $c_1,c_2,c_3\in A$. For example, tip points labeled in figure~\ref{fig:1} get reduced to the following algebraic expressions:
\begin{align*}
	\phi(\overline{1})=\phi(1111\dots)=1+c_1+c_1^2+c_1^3+\dots=\sum_{k=0}^\infty c_1^k&=\frac{1}{1-c_1},\\
	\phi(23\overline{1})=1+c_2+c_2c_3+c_2c_3c_1+c_2c_3c_1^2+\dots=1+c_2+c_2c_3\sum_{k=0}^\infty c_1^k&=1+c_2+\frac{c_2c_3}{1-c_1},\\
	\phi(122\overline{1})=1+c_1+c_1c_2+\frac{c_1c_2^2}{1-c_1},\quad \phi(32\overline{1})=1+c_3+\frac{c_3c_2}{1-c_1}, \quad \phi(133\overline{1})&=1+c_1+c_1c_3+\frac{c_1c_3^2}{1-c_1}.
	\end{align*}	

\noindent The set of all tip points $\phi(w)$ is called the \textbf{tipset}~$F_A$ of the complex tree~$T_A$. A tipset~$F_A$ is a self-similar set
\begin{equation}
	F_A=f_1(F_A)\cup f_2(F_A)\cup f_3(F_A)
\end{equation}
generated by an iterated function system composed of~three contractive mappings $f_1$,~$f_2$,~and~$f_3$ defined as
\begin{equation}
	f_1(z)=1+c_1z \text{ , }\quad f_2(z)=1+c_2z \ \text{, and }\ f_3(z)=1+c_3z
\end{equation}
where $\{c_1,c_2,c_3\}=A$ and $0<|c_1|,|c_2|,|c_3|<1$. 
\noindent The self-similar nature of the tipset~$F_A$ implies that the tip-to-tip intersections between the three first-level pieces $f_1(F_A)$,~$f_2(F_A)$,~and~$f_3(F_A)$ is the only piece of information needed for capturing the topological structure of a tipset~$F_A$. We define the \textbf{topological set}~$Q_A$ of a complex tree $T_A$ as the following set of tip-to-tip equivalence relations
\begin{equation}
	Q_A:=\{a\sim b : \phi(a)=\phi(b)=x\in f_j(F_A)\cap f_k(F_A) \ \text{ where } j=a_1\neq b_1=k\}.
\end{equation}
 For example, the topological set of ternary trees in figures~\ref{fig:0}-\ref{fig:1} is $Q_A=\emptyset$ and $Q_A=\{23\overline{1}\sim 122\overline{1},32\overline{1}\sim 133\overline{1}\}$ respectively since the tipset of the first one is topologically homeomorphic to a Cantor set, and for the second one, the intersection of first-level pieces $f_1(F_A)\cap f_2(F_A)=\{x\}$~and~$f_1(F_A)\cap f_3(F_A)=\{x'\}$ is a pair of singletons uniquely encoded by $x=\phi(23\overline{1})=\phi(122\overline{1})$~and~$x'=\phi(32\overline{1})=\phi( 133\overline{1})$.

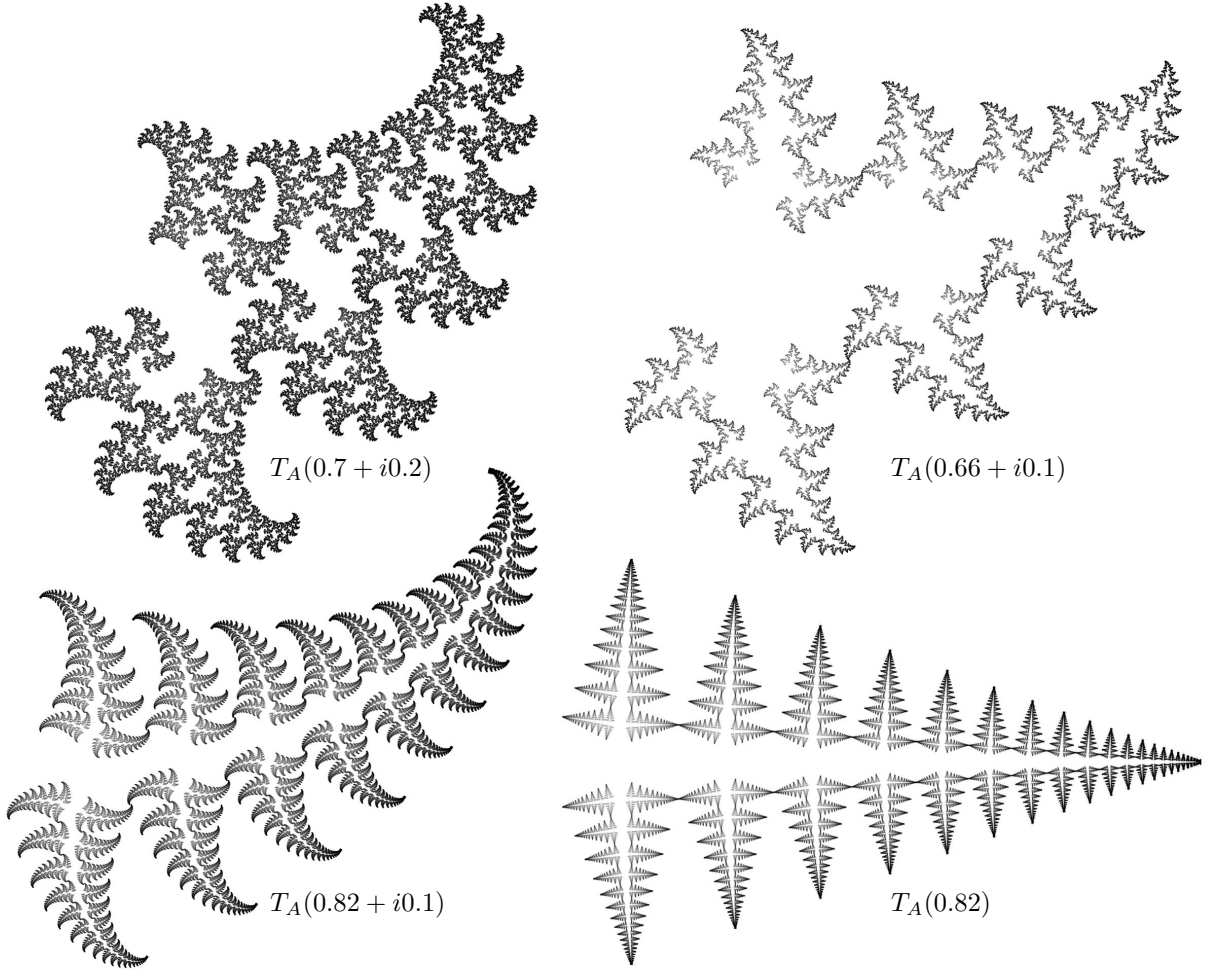
\begin{figure}[h!tbp]
	\centering
 \begin{overpic}[width=1.\textwidth,tics=10]{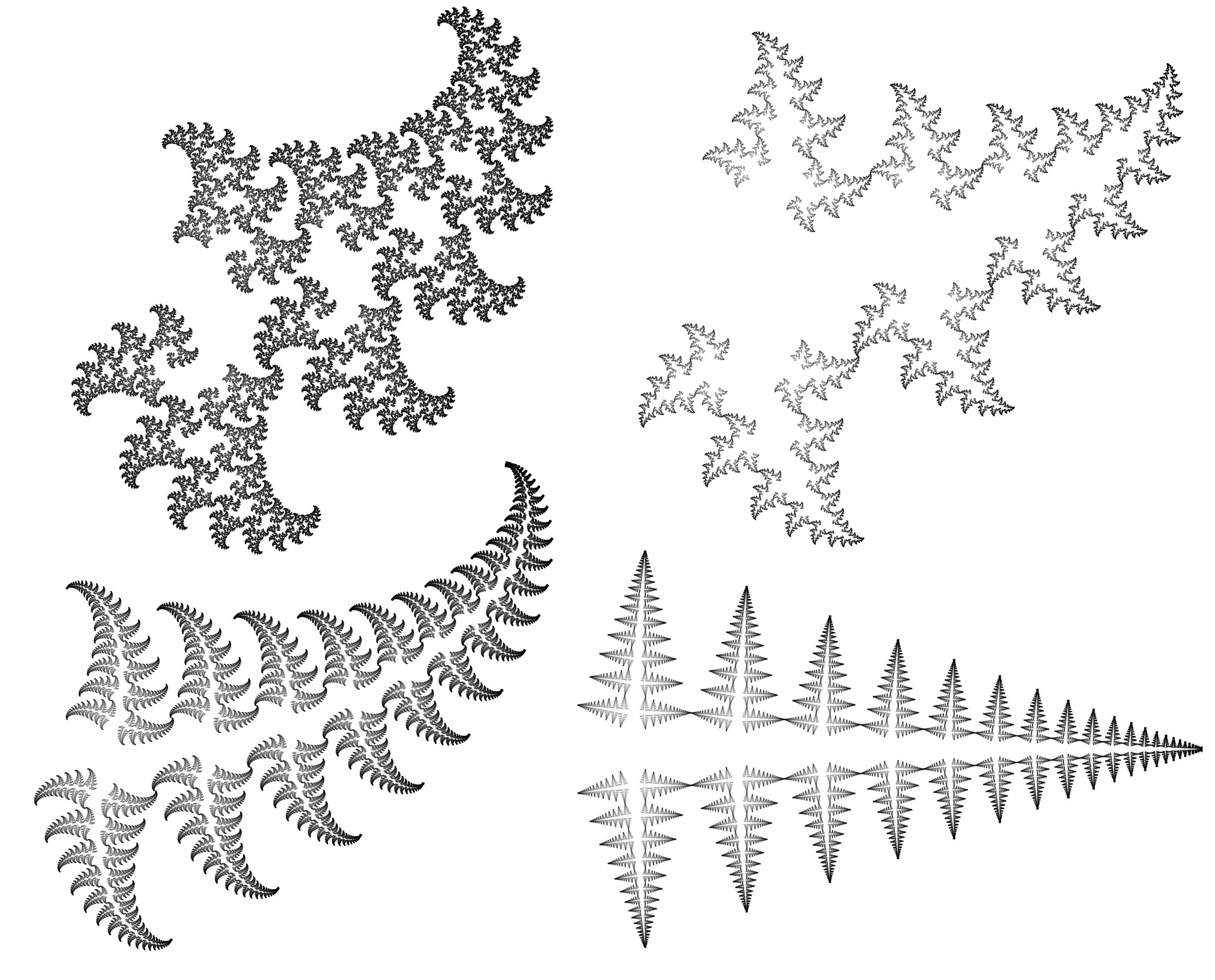}
 \put (24,42) {$T_A(0.7 + i0.2)$}
 \put (74,42) {$T_A(0.66 + i0.1)$}
 \put (24,7) {$T_A(0.82 + i0.1)$}
 \put (74,7) {$T_A(0.82)$}
\end{overpic}
	\caption{Topologically homeomorphic tipsets of stable trees $T_A(z)$ from the family obtained in~(\ref{familyeq}).}
	\label{fig:2}
\end{figure}

\section*{A Family of Fern-Like Connected Self-Similar Sets}
\noindent The method introduced in~\cite{Espigule2019} to obtain families of connected self-similar sets from topological sets~$Q_A$ of certain complex trees~$T_A$ can be applied for the tree depicted in figure~\ref{fig:1}. If we consider the topological set $Q_A=\{23\overline{1}\sim 122\overline{1},32\overline{1}\sim 133\overline{1}\}$ as a system of two equations $\{\phi(23\overline{1})=\phi(122\overline{1}),\phi(32\overline{1})=\phi( 133\overline{1})\}$ with letters $c_1,c_2,c_3\in A$ set as three unknown complex variables, then the system admits a parametric solution in one complex variable with~$c_1=z$ as unknown:
\begin{equation}\label{familyeq}
	\begin{cases}
	\phi(23\overline{1})=\phi(122\overline{1})\\
	\phi(32\overline{1})=\phi(133\overline{1})
	\end{cases}=\begin{cases}
	c_2+\frac{c_2c_3}{1-z}=z+zc_2+\frac{zc_2^2}{1-z}\\
	c_3+\frac{c_3c_2}{1-z}=z+zc_3+\frac{zc_3^2}{1-z}
	\end{cases}
	\  
	\begin{cases}
	c_2(z):=\frac{1-z-z^2+z^3+\sqrt{1-2 z-z^2-z^4+2 z^5+z^6}}{2 (z+z^2)}\\
	c_3(z):=2+(1+z^2)/z-c_2(z)\\
	\end{cases}
\end{equation}
\noindent This one-parameter family $T_A(z)=T\{z,c_2(z),c_3(z)\}$ is defined for~$\mathcal{R}:=\{z\in\mathbb{C}:0<|z|,|c_2(z)|,|c_3(z)|<1\}$. Notice that by construction we have that $\{23\overline{1}\sim 122\overline{1},32\overline{1}\sim 133\overline{1}\}\subseteq Q_A(z)$ for all $z\in\mathcal{R}$, i.e. all tipsets $F_A(z)$ are connected. Therefore the connectivity locus for this family of self-similar sets is the entire region~$\mathcal{R}$.

\begin{figure}[h!tbp]
	\centering
	\includegraphics[width=6.2in]{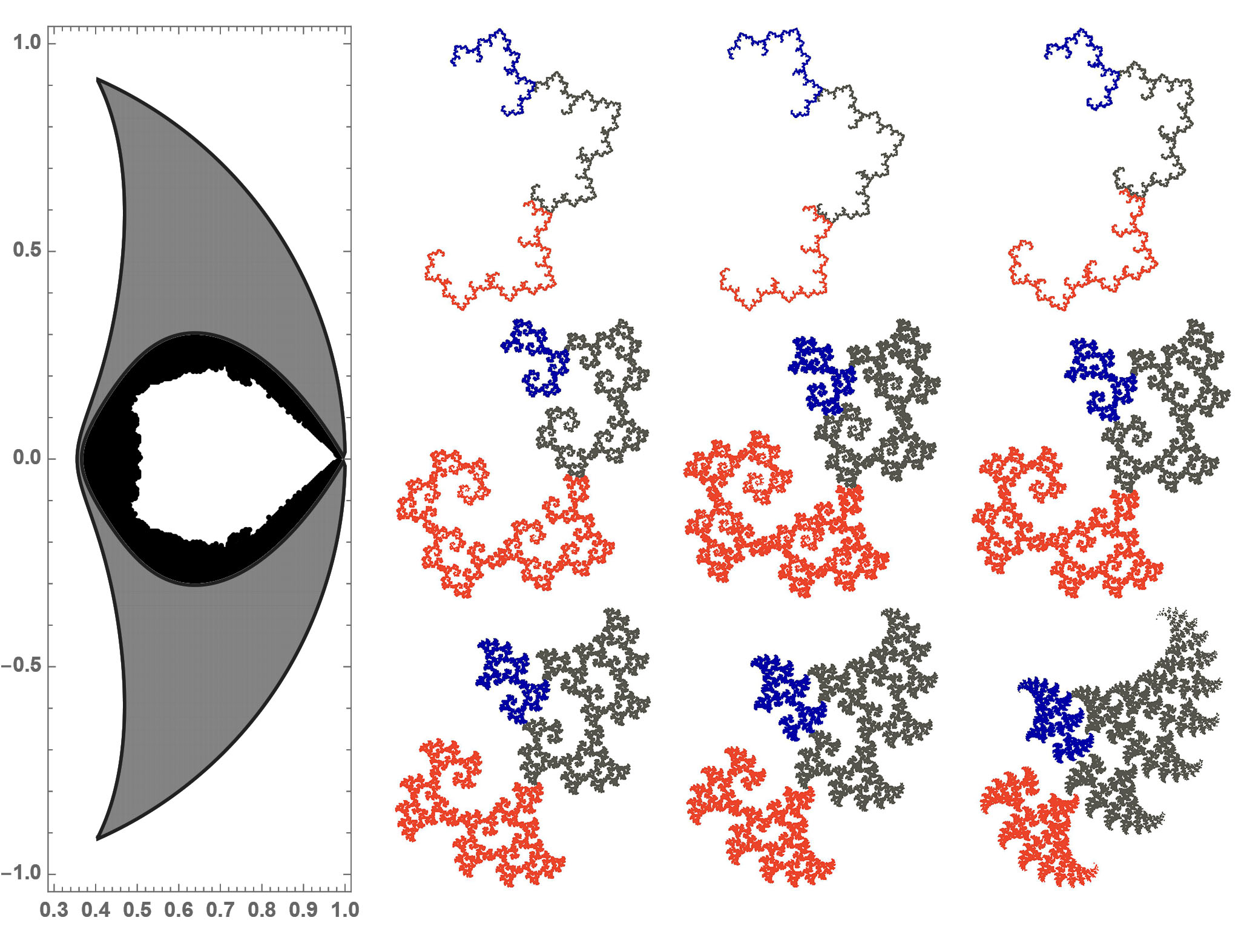}
	\caption{The unstable set $\mathcal{M}\subset\mathcal{R}$ with six tipsets of unstable trees $T_A(z)$ found in the boundary $\partial\mathcal{M}$.}
	\label{M}
\end{figure}

\section*{The Unstable Set~$\mathcal{M}$}

\begin{figure}[h!tbp]
	\centering
	\includegraphics[width=6.6in]{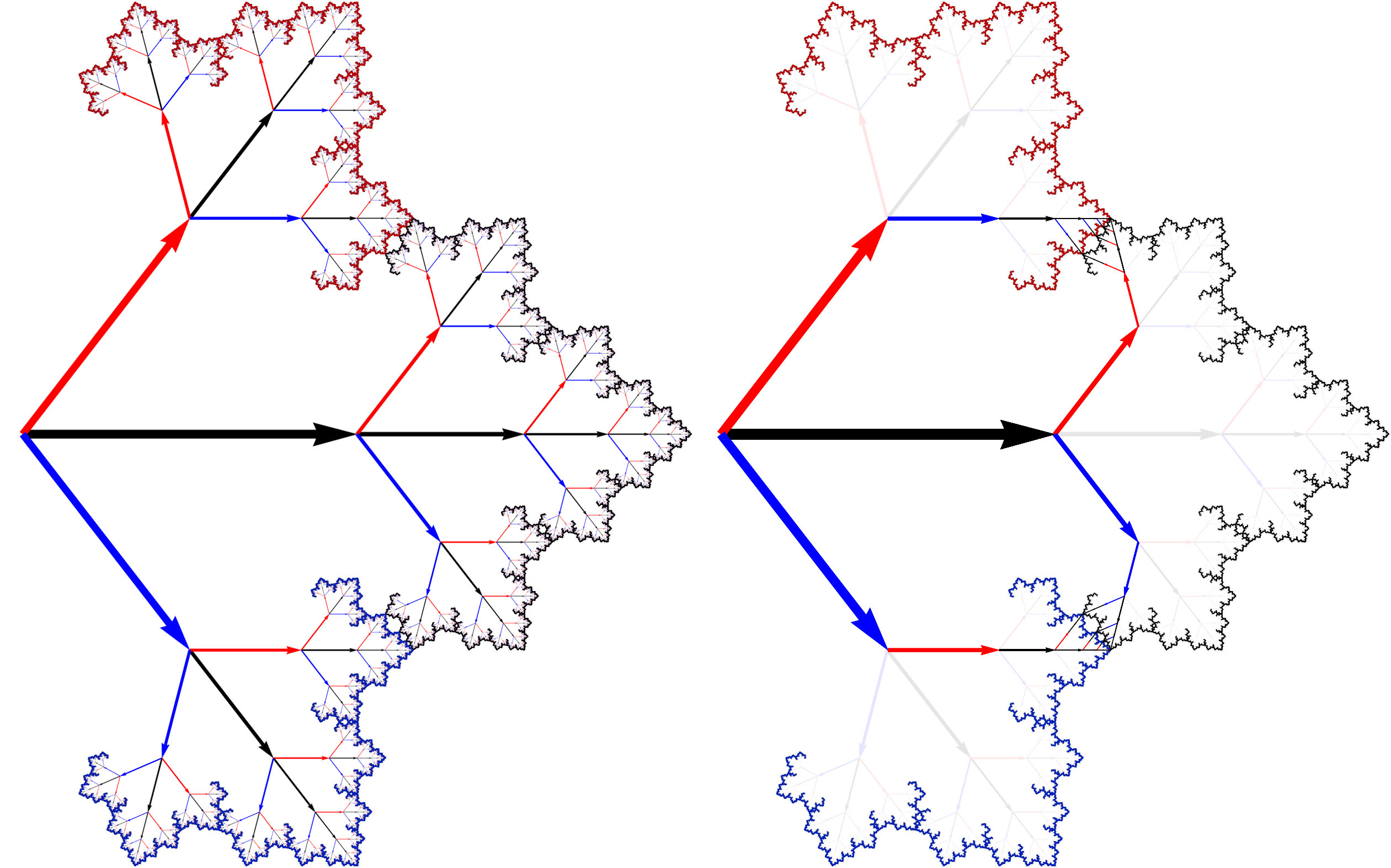}
	\caption{Unstable tree $T_A(1/2)=T\{1/2,1/4+i\sqrt{15}/12,1/4-i\sqrt{15}/12\}$ which is also mirror-symmetric.}
	\label{fig:3}
\end{figure}

A much more informative map about the topological structure of the family is given by what we call the unstable set defined as $\mathcal{M}:=\{z\in\mathcal{R}: Q_A(z)\neq\{23\overline{1}\sim 122\overline{1},32\overline{1}\sim 133\overline{1}\} \}$. The complement of~$\mathcal{M}$ is the stable set~$\mathcal{K}=\mathcal{R}\backslash\mathcal{M}$ entirely composed of points~$z$ with tipsets $F_A(z)$ topologically homeomorphic to $F_A(1/\tau)=F_A(-1/2+\sqrt{5}/2)$ shown in figure~\ref{fig:1}, i.e. $Q_A(z)=\{23\overline{1}\sim 122\overline{1},32\overline{1}\sim 133\overline{1}\}$. Such tipsets are structurally stable because we can always find an~$\epsilon$-neighborhood of $z\in\mathcal{K}$ with constant topological set $Q_A(z+\epsilon)=Q_A(z)=\{23\overline{1}\sim 122\overline{1},32\overline{1}\sim 133\overline{1}\}$.
 On the other hand, for $z\in\mathcal{M}$, such neighborhoods do not exist and the tipsets $F_A(z)$ are structurally unstable, any perturbation away from $z$ will destroy the original topological set $Q_A(z)\neq Q_A(z+\epsilon)$. Examples of stable and unstable trees~$T_A(z)$ are shown in figure~\ref{fig:2}~and~figures~\ref{M}-\ref{fig:3} respectively.

\noindent A direct way to determine points~$z$ in the unstable set~$\mathcal{M}$ with the aid of a computer software like \textit{Mathematica} consists in imposing an extra tip-to-tip equivalence relation $a\sim b\notin\{23\overline{1}\sim 122\overline{1},32\overline{1}\sim 133\overline{1}\}$ and then solving the equality $\phi(a)=\phi(b)$ with complex-valued letters $c_1$, $c_2$, and~$c_3$ replaced by those of the parametric alphabet $A(z)=\{1\rightarrow z,2\rightarrow c_2(z),3\rightarrow c_3(z)\}$. For example, with $2313\overline{1}\sim 1222\overline{1}$ we get~$z=1/2$ since $\phi(2313\overline{1})=\phi(1222\overline{1})$ gets reduced into:
\begin{align*}
1+c_2(z)+c_2(z)c_3(z)(1+z+zc_3(z)/(1-z))=1+z+zc_2(z)+zc_2(z)^2+zc_2(z)^2/(1-z)\\
0=(2 z-1)\left(-1+2 z+2 z^3+z^4+(z-1) \sqrt{1-2 z-z^2-z^4+2 z^5+z^6}\right)/(z+z^2).
\end{align*}
Therefore the mirror-symmetric tree $T_A(1/2)=T\{1/2,1/4+i\sqrt{15}/12,1/4-i\sqrt{15}/12\}$ shown in figure~\ref{fig:3} is unstable. Notice that the topological set of $T_A(1/2)$ has a numerable infinite number of tip-to-tip equivalence relations, $Q_A(1/2)=\{23\overline{1}\sim 122\overline{1},32\overline{1}\sim 133\overline{1},2313\overline{1}\sim 1222\overline{1},3212\overline{1}\sim 1333\overline{1},\dots\}$, but only one extra equivalence relation was actually needed to compute $z_0=1/2$. By automating this process of computing points $z_0\in\mathcal{M}$ from extra equivalence relations $a\sim b$ with $a_1\neq b_1$ we can approximate the unstable set~$\mathcal{M}$, see figure~\ref{M}. 

\begin{figure}[h!tbp]
	\centering
	\includegraphics[width=7in]{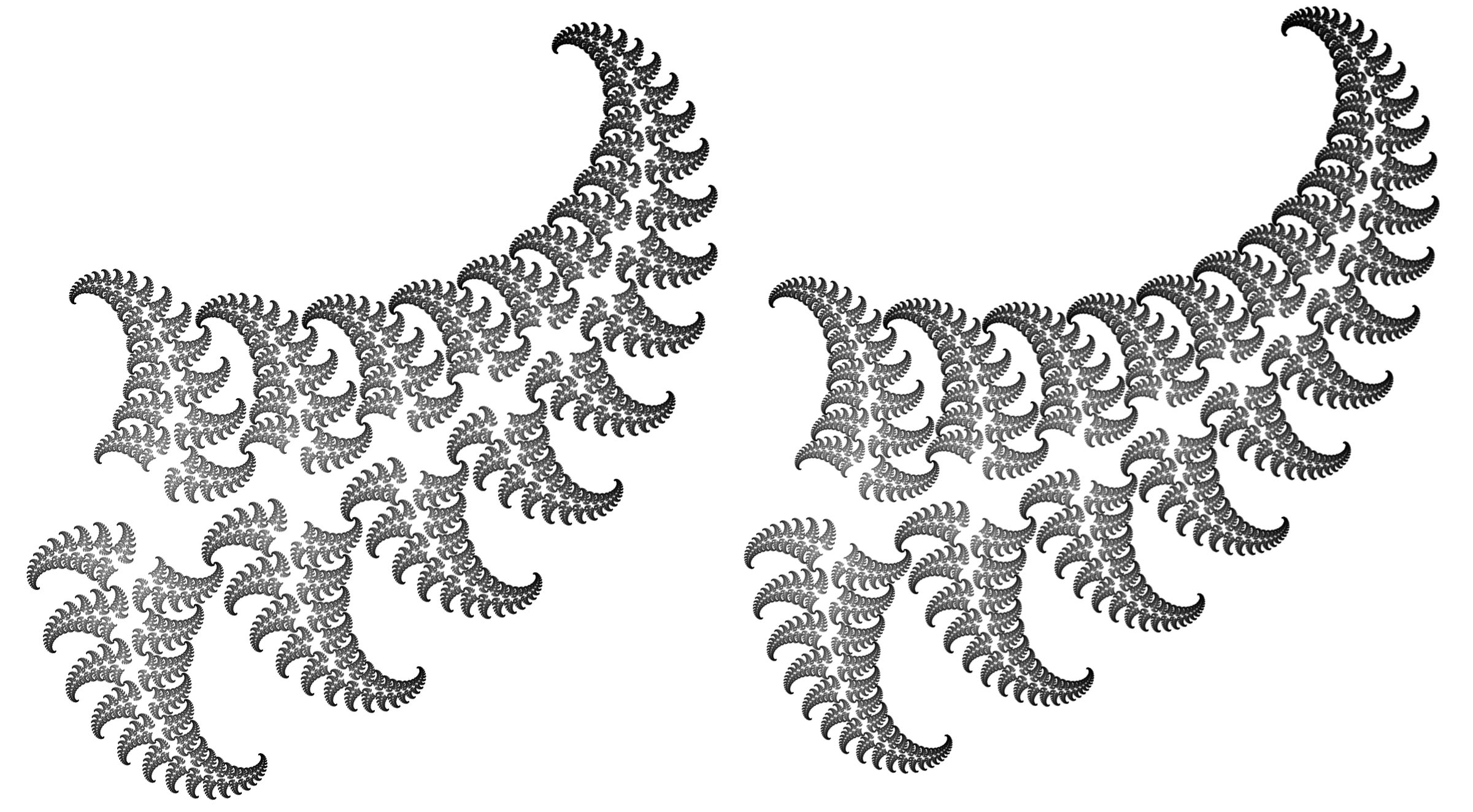}
	\caption{Tipsets $F_A(z)$~and~$F_A(z')$ where $z\approx0.814 + i0.146$ is a root of $1 - z^2 - z^3 - z^4+z^5 + z^6 + z^7$ and~$z'\approx 0.843 + i0.126$ is a root of $1 - z^2 - z^4 - z^5 + z^6 + z^7 + z^9$. Both points lie in the boundary~$\partial\mathcal{M}$ with extra tip-to-tip equivalence relations $21113\overline{1}\sim 12\overline{1}$ and~$211113\overline{1}\sim 12\overline{1}$.}
	\label{landmark}
\end{figure}

\noindent This brute-force method to compute $\mathcal{M}$ is exact but rather slow. A big chunk of the unstable set $\mathcal{M}$ is better obtained by an analytic method considered in~\cite{Espigule2019} which is based on the Hausdorff dimension and the open set condition. For our family of study we have that the analytic region $\mathcal{M}_2\subset \mathcal{M}$ is given by
\begin{equation}
\mathcal{M}_2:=\{z\in \mathcal{R} : 1<|z|^2+|c_2(z)|^2+|c_3(z)|^2\},
\end{equation}
see the gray region~$\mathcal{M}_2$ in figure~\ref{M}. For $z\notin\mathcal{M}$ the self-similar tipset~$F_A(z)$ satisfies the open set condition and their Hausdorff dimension coincides with the similarity dimension $\alpha$ which is the positive number satisfying
\begin{equation}\label{sim-value}
	|z|^\alpha+|c_2(z)|^\alpha+|c_3(z)|^\alpha=1.
\end{equation}
If $z\in\mathcal{M}$, then the Hausdorff dimension of a tipset $F_A(z)$ is not easy to compute in general and~eq.~(\ref{sim-value}) does not apply except in rare occasions when there are no overlaps and the pieces just-touch. This is precisely what happens for $z\in\partial\mathcal{M}$ where the open set condition still applies. Another remarkable property of the boundary $\partial\mathcal{M}$ is that extreme points penetrating into the stable set turn out to be interesting algebraic numbers, see the pair of examples in figure~\ref{landmark} with Hausdorff dimension $\alpha\approx1.763$ and~$\alpha\approx1.769$ respectively.

 \begin{figure}[h!tbp]
	\centering
		\begin{overpic}[width=1\textwidth,tics=10]{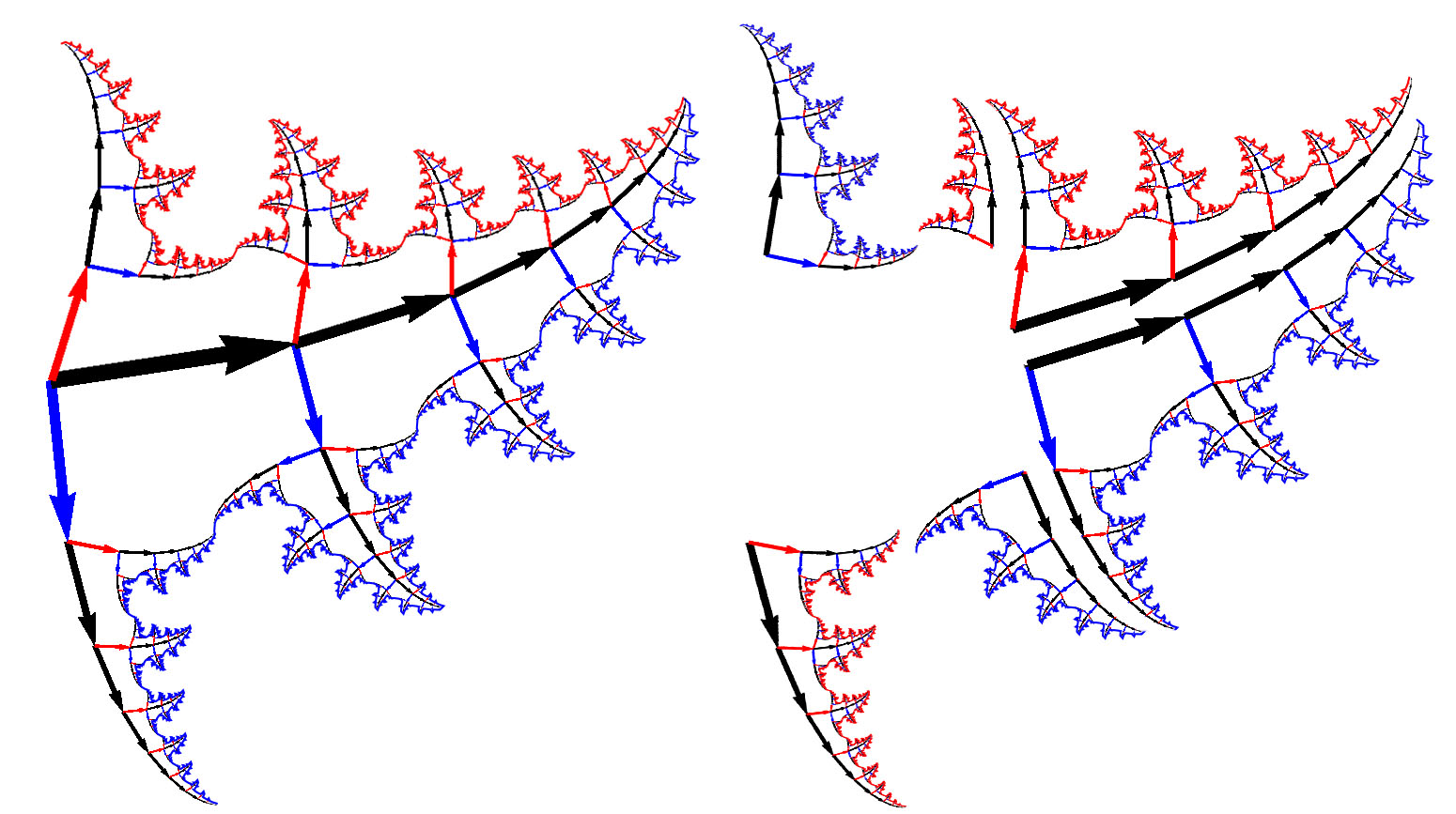}
 \put (24,49) {\huge$C$}
 \put (32,18) {\huge$D$}
 \put (0,54) {$\phi(2\overline{1})$}
 \put (45,50) {$\phi(\overline{1})$}
 \put (16,0) {$\phi(3\overline{1})$}
 \put (54,54) {\Large$f_2(D)$}
 \put (59,50) {\large$f_1(f_2(C))$}
 \put (74,50) {\Large$f_1(C)$}
 \put (82,18) {\Large$f_1(D)$}
 \put (68,9) {\large$f_1(f_3(D))$}
 \put (61,5) {\Large$f_3(C)$}
 \end{overpic}
	\caption{Shortest path $C\cup D\subset F_A(z)$ that goes from the tip point $\phi(2\overline{1})$ to the tip point $\phi(3\overline{1})$.}
	\label{pathAB}
\end{figure}
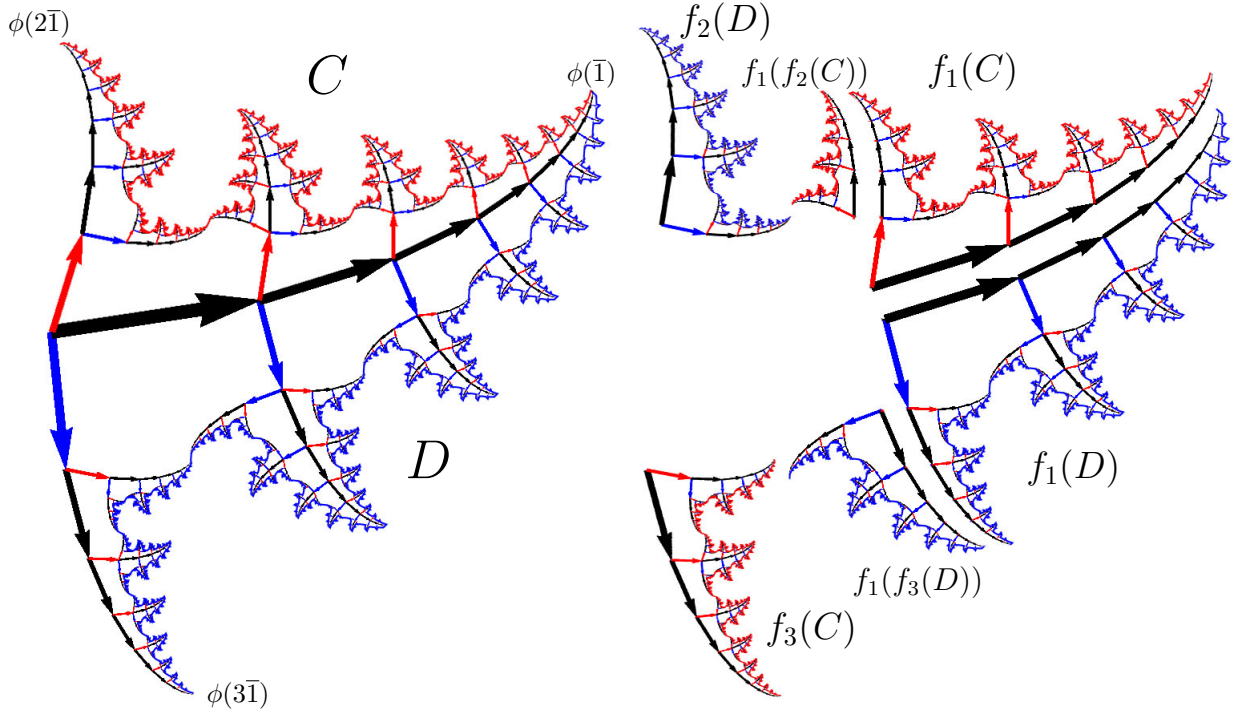

\section*{The Shortest Path from $\phi(2\overline{1})$ to $\phi(3\overline{1})$}

\begin{figure}[h!tbp]
	\centering
	\includegraphics[width=6.5in]{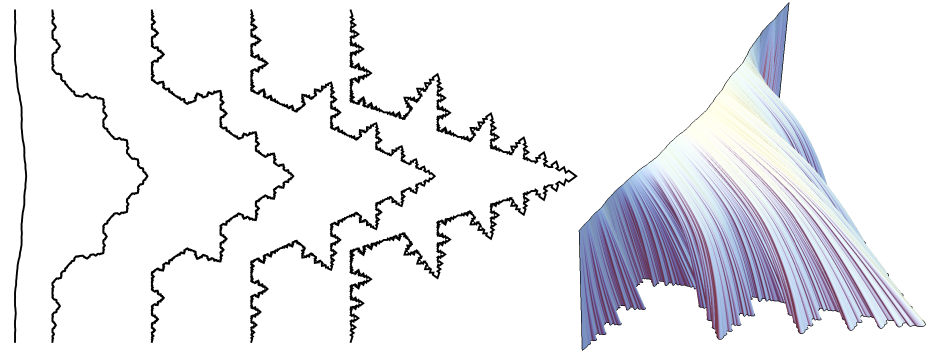}
	\caption{Shortest paths $C\cup D\subset F_A(z)$ for $z\in(1-\sqrt{2} ,1)$ with the 3D surface generated by them.}
	\label{surface}
\end{figure}

Christoph Bandt considered the notion of geodesics in self-similar sets, his method reported in~\cite{bandt2014geometry} can be adapted to tipsets~$F_A(z)$ for parameters $z$ in the stable set~$\mathcal{K}$ our family of study. The shortest path~$C\cup D\subset F_A(z)$ from the tip point $\phi(2\overline{1})$ to $\phi(3\overline{1})$ that goes through $\phi(\overline{1})\in C\cap D$ is given by $C= f_2(D) \cup f_1(f_2(C)) \cup f_1(C)$ and $D= f_1(D) \cup f_1(f_3(D))\cup f_3(C)$ where $C$ and $D$ are the curves depicted in figure~\ref{pathAB}. When we move the parameter~$z$ along the real line starting at $z=1-\sqrt{2}$, trees $T_A(z)$ are mirror-symmetric and the fractal dimension of $C\cup D$ starts to increase, see figure~\ref{surface}. As a final application, consider stacking these curves $C\cup D$ at heights given by the parameter $z\in(1-\sqrt{2},1)$, with $\phi(2\overline{1})$ and $\phi(3\overline{1})$ fixed at the same position. The result is the three-dimensional surface depicted in figure~\ref{surface}.

\section*{Summary and Conclusions}
The family considered in this paper represents a tiny sample of what is out there. The space of possible parametric families of tipset connected $n$-ary complex trees is incredibly vast and rich. Nonetheless we believe that the family covered here ranks hight in this space and it deserved to be considered apart. The theoretical basis set in~\cite{Espigule2019} provides a unified approach to previously known results on symmetric fractal trees and self-similar sets in general, see~\cite{bandt2014geometry}~\cite{Espigule2013} and references within.

\begin{figure}[h!tbp]
	\centering
	\includegraphics[width=7.in]{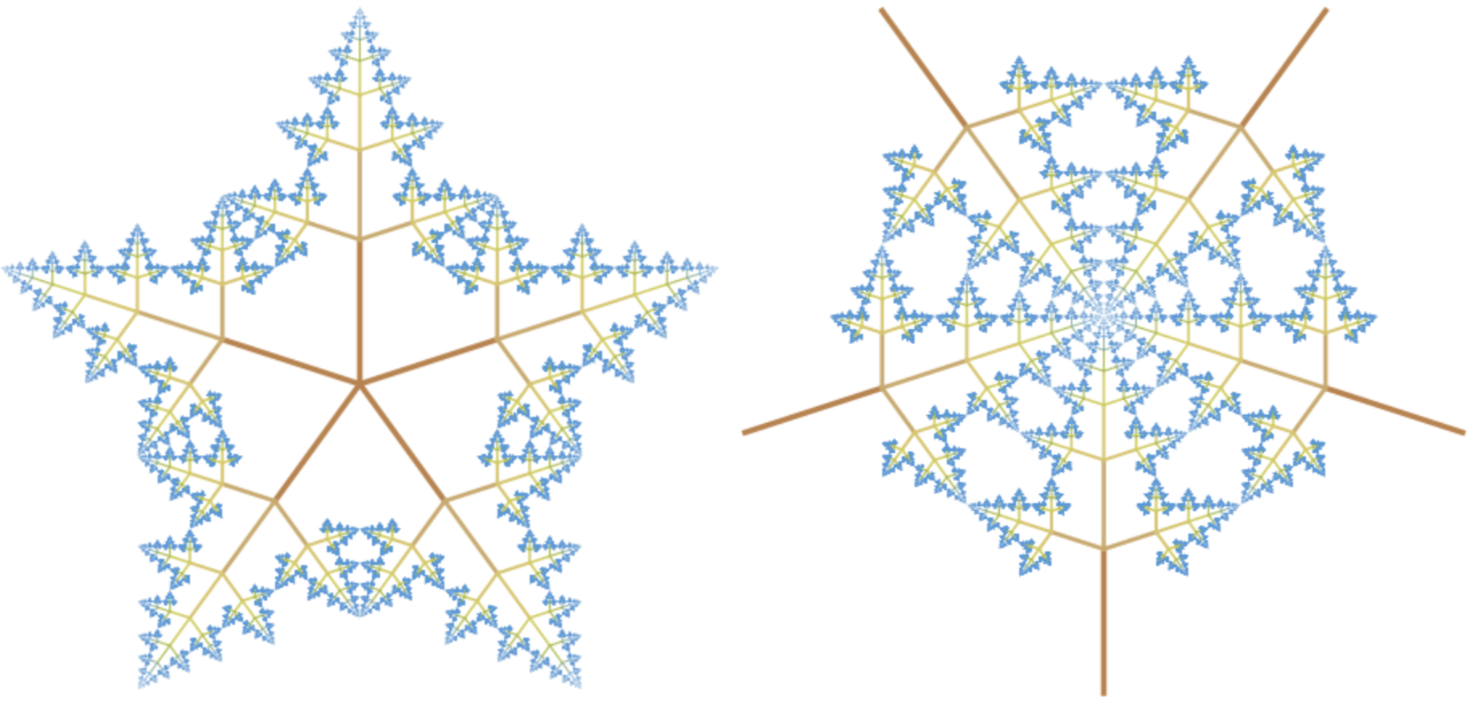}
	\caption{Five-fold rotational symmetry of the golden ternary tree~$T_A(1/\tau)$ shown in figure~\ref{fig:1}. Images reproduced from the author's post published in the Wolfram Blog~\cite{Espigule2014}.}
	\label{fivefold}
\end{figure}

\section*{Acknowledgements}
The author specially thanks the referees for suggesting several improvements, and IMUB's Holomorphic Dynamics group, N\'uria Fagella, Xavier Jarque, Toni Garijo, Robert Florido, etc for all their support, dedication, and valuable discussions during the author's research project on complex trees 2017-19. This project was partially supported by a research grant awarded by IMUB and the University of Barcelona. 
The author also thanks Wolfram Research and in particular Theodore Gray, Chris Carlson, Michael Trott, Vitaliy Kaurov, Todd Rowland, and Stephen Wolfram for showing interest on the author's early work and for being a source of inspiration when it comes to computational experiments~\cite{wolfram}. All the figures and diagrams of this paper were done with \textit{Mathematica}. Finally, the author would like to express his gratitude to Susanne Kr\"omker, Jofre Espigul\'e, Pere Pascual, Gaspar Orriols, Warren Dicks, Tara Taylor, Hans Walser, Robert Fathauer, Tom Verhoeff, Henry Segerman, Saul Schleimer, Michael Barnsley, Przemyslaw Prusinkiewicz, and many others for all their time and encouragement received from them during an early stage of the author's research~2012-14.

    
{\setlength{\baselineskip}{13pt} 
\raggedright				
\bibliographystyle{unsrt}

} 
   
\end{document}